\documentclass[12pt,dvipdfmx]{amsart}

\usepackage{amssymb,hyperref}

\newtheorem{theorem}{Theorem}[section]
\newtheorem{proposition}[theorem]{Proposition}
\newtheorem{lemma}[theorem]{Lemma}
\newtheorem{corollary}[theorem]{Corollary}

\theoremstyle{definition}
\newtheorem{definition}[theorem]{Definition}

\newtheorem{example}[theorem]{Example}

\newtheorem{remark}[theorem]{Remark}

\newcommand{\ZZ}{ \ensuremath{\mathbb{Z}}}

\newcommand{\lex}{{\mathrm{lex}}}

\newcommand{\Tor}{\ensuremath{\mathrm{Tor}}\hspace{1pt}}

\def\cocoa{{\hbox{\rm C\kern-.13em o\kern-.07em C\kern-.13em o\kern-.15em A}}}

\newcommand{\B}{\mathcal{B}}
\newcommand{\bier}{\mathrm{Bier}}
\newcommand{\cc}{\mathbf c}
\newcommand{\pol}{\mathrm{pol}}
\newcommand{\aaa}{\mathbf a}
\newcommand{\bb}{\mathbf b}
\newcommand{\facet}{\mathrm{Facet}_\cc}
\newcommand{\lk}{\mathrm{lk}}

\begin{document}

\title{Spheres arising from multicomplexes}

\author{Satoshi Murai}
\address{
Satoshi Murai,
Department of Mathematical Science,
Faculty of Science,
Yamaguchi University,
1677-1 Yoshida, Yamaguchi 753-8512, Japan
}



\begin{abstract}
In 1992, Thomas Bier introduced a surprisingly simple way to construct a large number of simplicial spheres.
He proved that, for any simplicial complex $\Delta$ on the vertex set $V$ with $\Delta \ne 2^V$, the deleted join of $\Delta$ with its Alexander dual $\Delta^\vee$ is a combinatorial sphere.
In this paper,
we extend Bier's construction to multicomplexes,
and study their combinatorial and algebraic properties.
We show that all these spheres are shellable and edge decomposable,
which yields a new class of many shellable edge decomposable spheres that are not realizable as polytopes.
It is also shown that
these spheres are related to polarizations and Alexander duality for monomial ideals which appear in commutative algebra theory.
\end{abstract}

\maketitle

\section*{Introduction}

In an unpublished work in 1992, Thomas Bier presented a surprisingly simple way to construct a large number of simplicial spheres
by using Alexander duality.
Let $\Delta$ be a simplicial complex on $[n]=\{1,2,\dots,n\}$ which is not the set of all subsets of $[n]$,
and let $\Delta^{\vee}=\{ F \subset [n]: [n] \setminus F \not \in \Delta \}$ be the Alexander dual of $\Delta$.
The \textit{Bier sphere $\bier (\Delta)$ of $\Delta$}
is defined as the deleted join of $\Delta$ and $\Delta^\vee$,
in other words, $\bier (\Delta)$ is the simplicial complex on $\{x_1,\dots,x_n,y_1,\dots,y_n\}$
defined by
$$
\bier (\Delta) =
\{ X_F \cup Y_G:
F \in \Delta,\
G \in \Delta^{\vee},\
F \cap G= \emptyset\},
$$
where $X_F=\{x_i:i \in F\}$ and $Y_G=\{y_i:i \in G\}$.
Bier proved that this simplicial complex is indeed a triangulation of a sphere,
and Bier spheres are of interest in topological combinatorics
in connection with the van Kampen-Flores theorem and non-polytopal triangulations of spheres \cite{BPSZ,Ma}.

The main purpose of this paper is to extend
Bier's construction to finite multicomplexes.
Our construction method is quite different from the original approach given by Bier.
Our approach is similar to the construction of Billera-Lee polytopes \cite{BL} and Kalai's squeezed spheres \cite{Ka}.
We first define a shellable ball $\B(M)$ from a finite multicomplex $M$,
and define the Bier sphere of $M$ as its boundary.

Here we briefly define our generalized Bier spheres.
Fix $\cc=(c_1,\dots,c_n) \in \ZZ^n_{\geq 0}$.
A monomial $x_1^{a_1}x_2^{a_2} \cdots x_n^{a_n}$, where $x_1,x_2,\dots,x_n$ are variables,
is called a \textit{$\cc$-monomial} if $a_i \leq c_i$ for all $i$.
A {\em $\cc$-multicomplex} $M$ is a non-empty set of $\cc$-monomials such that if 
$u \in M$ and a monomial $v$ divides $u$ then $v \in M$.
Let
$$\widetilde X_i=\{ x_i^{(0)},x_i^{(1)},\dots,x_i^{(c_i)}\}$$
be sets of new indeterminates for $i=1,2,\dots,n$ and let $\widetilde X=\bigcup_{i=1}^n \widetilde X_i$.
We define the simplicial complex $\B_\cc(M)$ to be the simplicial complex generated by
$$\left\{ \widetilde X \setminus \{ x_1^{(a_1)},x_2^{(a_2)},\dots,x_n^{(a_n)} \}: x_1^{a_1} x_2^{a_2}\cdots x_n^{a_n} \in M \right \}.$$
In Lemma \ref{ball},
we prove that if $M$ is not the set of all $\cc$-monomials
then $\B_\cc(M)$ is a $(|\cc|-1)$-dimensional shellable ball,
where $|\cc|=c_1+ \cdots +c_n$.
We define the
\textit{Bier sphere $\bier_\cc (M)$ of a multicomplex $M$ (w.r.t.\ $\cc$)}
as the boundary of $\B_\cc(M)$.

The above definition looks very different from the construction of classical Bier spheres of simplicial complexes.
However, we show that classical Bier spheres are Bier spheres of multicomplexes w.r.t.\ $\cc=(1,1,\dots,1)$.
Note that another generalization of Bier's construction (to posets) was given by Bj\"orner et al.\ \cite{BPSZ}.

We also
prove the following combinatorial properties of $\bier_\cc(M)$
(the definition of $g$-vectors, shellability and edge decomposability will be given in Sections 1 and 4).
\begin{itemize}
\item $\bier_\cc(M)$ is shellable (Theorem \ref{shellable}).
\item The $g$-vector of $\bier_\cc(M)$ is given by
$$
g_i(\bier_\cc(M))= \# \{ u \in M : \deg u =i\} -\# \{u \in M: \deg u=|\cc|-i\},
$$
where $\# X$ denotes the cardinality of a finite set $X$ (Theorem \ref{hvector}).
\item $\bier_\cc(M)$ is edge decomposable (Theorem \ref{edgedecomposable}).
\end{itemize}

The first and second results generalize
the results of Bj\"oner et al.\ \cite[Theorems 4.1 and 5.2]{BPSZ}
who proved the same statement for classical Bier spheres of simplicial complexes.
Edge decomposability was introduced by Nevo \cite{Ne} in the study of $g$-vectors of simplicial spheres.
This property is important since if a simplicial complex is edge decomposable
then its face vector satisfies McMullen's $g$-condition.
Thus the result yields a new class of simplicial spheres whose face vectors satisfy
McMullen's $g$-condition,
since most Bier spheres are not realizable as polytopes \cite{BPSZ,Ma}.
It would be of interest to find a simple proof of this fact by using the second result
without using edge decomposability.

The construction of $\bier(M)$ is inspired by the study of shellability of multicomplexes and polarization
which are techniques in commutative algebra theory \cite{HP,So}.
Let $S=K[x_1,\dots,x_n]$ be a polynomial ring over a field $K$ with $\deg x_i=1$ for any $i$.
Let $\overline \cc =\cc + (1,1,\dots,1)$.
For a $\overline \cc$-monomial $x_1^{a_1}x_2^{a_2}\cdots x_n^{a_n}$,
its {\em polarization (w.r.t.\ $\overline \cc$)} is the squarefree monomial
$$\pol(x_1^{a_1}x_2^{a_2} \cdots x_n^{a_n})=
\prod_{a_i \ne 0} (x_{i,0} x_{i,1}\cdots x_{i,a_i-1})$$
in a polynomial ring $\widetilde S=K[x_{i,j}: 1 \leq i \leq n,\ 0 \leq j \leq c_i]$.
For a monomial ideal $I \subset S$ minimally generated by $\overline \cc$-monomials $u_1,\dots,u_s$,
the {\em polarization of $I$ (w.r.t.\ $\overline \cc$)} is the squarefree monomial ideal
$$\pol(I)=\big(\pol(u_1),\dots,\pol(u_s)\big) \subset \widetilde S.$$
Polarizations of monomial ideals are used in the study of graded Betti numbers in commutative algebra.
An advantage of polarization is that,
since taking polarizations does not change graded Betti numbers,
polarization sometimes reduce an algebraic problem on graded Betti numbers of monomial ideals
to a combinatorial problem of simplicial complexes.
Our key observation which connects polarization and Bier spheres is the following.
Let $M$ be a $\cc$-multicomplex,
and let $I(M)$ be the ideal generated by all monomials in $S$ which are not in $M$.
In Lemma \ref{jahan},
we show
$$
\pol\big(I(M)\big)=\left(\prod_{x_i^{(j)} \in F} x_{i,j}: F \subset \widetilde X,\ F \not \in \B_\cc(M)\right).$$ 
Thus we show that the Stanley-Reisner ideal of the Bier ball $\B_\cc(M)$
is the polarization of the monomial ideal $I(M)$.

By using the above fact,
we study algebraic aspects of Bier spheres.
Since Bier spheres of simplicial complexes are defined by using Alexander duality,
it is natural to ask if Bier spheres of multicomplexes are also related to Alexander duality.
In Section 3, we show that Bier spheres of multicomplexes can be defined by using
Alexander duality for monomial ideals, introduced by Miller \cite{Mi,Mi2}.
In Section 5, we discuss a connection between Bier spheres and linkage theory,
and study graded Betti numbers of Stanley-Reisner rings of Bier spheres.

The results of the paper show that Bier spheres and Kalai's squeezed spheres \cite{Ka} can be constructed in a similar way.
Kalai's squeezed spheres also arise from finite multicomplexes by certain operations
(see \cite[p.\ 6]{Ka} and \cite[Proposition 4.1]{Mu1}),
and give many shellable edge decomposable spheres which are not realizable as polytopes (\cite{Ka,Le,Mu2}).
It might be of interest to find a general construction of shellable spheres which includes both Bier spheres and Kalai's squeezed spheres.

\section{Bier sphere of a multicomplex}

We first recall the basics on simplicial complexes and multicomplexes.

Let $V$ be a finite set.
A simplicial complex $\Delta$ on $V$ is a collection of subsets of $V$ such that
if $F \in \Delta$ and $G \subset F$ then $G \in \Delta$
(we do not assume that $\Delta$ contains all $1$-subsets of $V$).
An element $F \in \Delta$ with $\#F=i+1$ is called an {\em ($i$-dimensional) face of $\Delta$},
and maximal faces under inclusion are called {\em facets}.
The \textit{dimension of $\Delta$} is the maximal dimension of its faces.
A simplicial complex is said to be {\em pure} if all its facets have the same cardinality.
For subsets $F_1,F_2,\dots,F_s$ of $V$,
we write $\langle F_1,F_2,\dots,F_s\rangle$ for the simplicial complex generated by $F_1,F_2,\dots,F_s$,
in other words,
$$\langle F_1,F_2,\dots,F_s\rangle=\big\{ G \subset V: G \subset F_i \mbox{ for some }i \in \{1,2,\dots,s\}\big\}.$$

\begin{definition}
A $(d-1)$-dimensional pure simplicial complex $\Delta$ is said to be {\em shellable} if there is an order $F_1,F_2,\dots,F_s$
of the facets of $\Delta$ such that
$$\langle F_1,F_2,\dots,F_{i-1} \rangle \cap \langle F_i \rangle$$
is generated by subsets of $F_i$ of cardinality $d-1$.
The order $F_1,F_2,\dots,F_s$ is called a {\em shelling of $\Delta$}.
\end{definition}

We say that a simplicial complex $\Delta$ is a {\em simplicial $d$-ball} (or {\em $d$-sphere})
if its geometric realization is homeomorphic to a $d$-ball (or $d$-sphere).
It is well-known that if a $(d-1)$-dimensional simplicial complex $\Delta$ is shellable and if any $(d-2)$-dimensional face of $\Delta$ is contained in at most two facets, then $\Delta$ is a simplicial ball or sphere (see \cite[Theorem 11.4]{Bj}).
This implies the following fact.

\begin{lemma}
\label{shellball}
Let $\Delta$ be a simplicial $(d-1)$-sphere and let $\Gamma$ be a $(d-1)$-dimensional subcomplex of $\Delta$ with $\Gamma \ne \Delta$.
If $\Gamma$ is shellable
then $\Gamma$ is a simplicial $(d-1)$-ball.
\end{lemma}

Let $X=\{x_1,x_2,\dots,x_n\}$ be a set of indeterminates.
For $\aaa=(a_1,a_2,\dots,a_n)\in \ZZ_{\geq 0}^n$,
we write
$$x^\aaa = x_1^{a_1}x_2^{a_2}\cdots x_n^{a_n}.$$
Fix $\cc=(c_1,c_2,\dots,c_n) \in \ZZ_{\geq 0}^n$.
A $\cc$-multicomplex is said to be \textit{$\cc$-full}
if it is the set of all $\cc$-monomials.
A $\cc$-multicomplex which is not $\cc$-full is called a \textit{proper $\cc$-multicomplex}.
Let
$$\widetilde X_i=\{ x_i^{(0)},x_i^{(1)},\dots,x_i^{(c_i)}\}$$
be sets of new indeterminates for $i=1,2,\dots,n$, and let
$$\widetilde X= \widetilde X_1 \bigcup \widetilde X_2 \bigcup \cdots \bigcup \widetilde X_n.$$

\begin{definition}
For any $\cc$-monomial $x^\aaa=x_1^{a_1}x_2^{a_2}\cdots x_n^{a_n}$,
let
$$F_\cc(x^\aaa)=
\widetilde X \setminus \{ x_1^{(a_1)},x_2^{(a_2)},\dots,x_n^{(a_n)} \}.
$$
For any $\cc$-multicomplex $M$,
we define the simplicial complex $\B_\cc(M)$ on $\widetilde X$ by
$$\B_\cc(M)= \left\langle F_\cc(x^\aaa): x^\aaa \in M \right \rangle.$$
\end{definition}

For simplicial complexes $\Delta$ and $\Gamma$ with disjoint vertices,
the simplicial complex
$$\Delta * \Gamma =\{F \cup G: F \in \Delta \mbox{ and } G \in \Gamma\}$$
is called the {\em join of $\Delta$ and $\Gamma$}.
Let $\partial \widetilde X_i=\{F \subset \widetilde X_i: F \ne \widetilde X_i\}$ for $i=1,2,\dots,n$,
and let
$$\Lambda_\cc=\partial \widetilde X_1 * \partial \widetilde X_2 * \cdots * \partial \widetilde X_n.$$
Then $\Lambda_\cc$ is a simplicial $(|\cc|-1)$-sphere,
where $|\cc|=c_1+\cdots+c_n$,
since each $\partial \widetilde X_i$ is the boundary of a simplex and since the join of two simplicial spheres is again a simplicial sphere.
Since $F_\cc(x^\aaa)=\widetilde X \setminus \{x_1^{(a_1)},\dots,x_n^{(a_n)}\}$ is a facet of $\Lambda_\cc$,
$\B_\cc(M)$ is a subcomplex of $\Lambda_\cc$ having the same dimension as $\Lambda_\cc$.

\begin{lemma} \label{ball}
Let $M$ be a $\cc$-multicomplex.
\begin{itemize}
\item[(i)] $\mathcal B_\cc(M) = \Lambda_\cc$ if and only if $M$ is $\cc$-full.
\item[(ii)] $\B_\cc(M)$ is shellable.
\item[(iii)] If $M$ is not $\cc$-full then $\B_\cc(M)$ is a simplicial $(|\cc|-1)$-ball.
\end{itemize}
\end{lemma}

\begin{proof}
(i) is straightforward.

(ii) was essentially proved in \cite[Theorem 4.3]{So}.
But we include a proof for the sake of completeness.
Let $x^\aaa \in M$ be a monomial that does not divide any other monomial in $M$.
It is enough to prove that
\begin{align}
\label{1.1}
\B_\cc(M\setminus \{x^\aaa\})\cap \langle F_\cc(x^\aaa) \rangle
=\big \langle F_\cc(x^\aaa) \setminus \{x_i^{(j)}\}: i=1,2,\dots,n,\ 0 \leq  j \leq a_i-1 \big\rangle.
\end{align}

We first show that the left-hand side contains the right-hand side.
Observe that $x^\aaa x_i^{-a_i+j} \in M$ if $j \leq a_i -1$.
Then $F_\cc(x^\aaa) \setminus \{x_i^{(j)}\} \subset F_\cc (x^\aaa x_i^{-a_i+j}) \in \B_\cc(M\setminus \{x^\aaa\})$.
Second, we show that the right-hand side contains the left-hand side.
It is enough to prove that
$$G=\{ x_i^{(j)}: i=1,2,\dots,n,\ 0\leq j \leq a_i-1\}$$
is not contained in $\B_\cc(M\setminus\{x^\aaa\})$
since the set $G$ is the (unique) smallest element among the elements in $\langle F_\cc(x^\aaa)\rangle$ which are not contained in the right-hand side of \eqref{1.1}.
Suppose contrary that $G \in \B_\cc(M\setminus\{x^\aaa\})$.
There exists a monomial $x^\bb \in M$ with $x^\bb \ne x^\aaa$ such that $G \subset F_\cc(x^\bb)$.
Then, by the definition of $F_\cc(-)$, we have $b_i \geq a_i$ for all $i$,
which implies $x^\aaa$ divides $x^\bb$.
This contradicts the choice of $x^\aaa$.

(iii)
By (i),
$\B_\cc(M) \subsetneq \Lambda_\cc$.
Then Lemma \ref{shellball} and (ii) say that $\B_\cc(M)$ is a simplicial $(|\cc|-1)$-ball.
\end{proof}

We call $\B_\cc(M)$ the {\em Bier ball of $M$ with respect to $\cc$.}

\begin{remark}\label{r1}
$\B_\cc(M)$ depends not only on $M$ but also on $\cc$.
However, the crucial case will be the case when $x_1^{c_1} \cdots x_n^{c_n}$
is equal to the least common multiple of monomials in $M$
since $\B_{(c_1+1,c_2,\dots,c_n)}(M)=\{x_1^{(c_1+1)}\}*\B_\cc(M)$ is a cone of $\B_\cc(M)$.
\end{remark}

\begin{remark}
Clearly $\B_{(c_1,\dots,c_{n-1},0)}(M)=\B_{(c_1,\dots,c_{n-1})}(M)$.
Thus one can assume $\cc \in \ZZ_{\geq 1}^n$.
However, in this paper we include $0$
since considering $\ZZ_{\geq 0}^n$ is convenient for induction purposes.
\end{remark}

\begin{remark}
In the special case when $M$ is the set of all monomials of degree $\leq k$ for some integer $k\geq 0$,
Lemma \ref{ball}(iii) was proved in \cite[Theorem 3.1]{HP}.
\end{remark}

\begin{remark}
Although we only prove shellability,
$\B_\cc(M)$ is vertex decomposable (see \cite[p.\ 1854]{Bj}).
Indeed, both the link and the deletion of $\B_\cc(M)$ w.r.t.\ $x_1^{(c_1)}$
are Bier balls.
\end{remark}

Now we define Bier spheres of multicomplexes.
Let $\Delta$ be a simplicial $(d-1)$-ball.
Then each $(d-2)$-dimensional face
of $\Delta$ is contained in at most two facets of $\Delta$.
Then its boundary
$$\partial \Delta = \langle F \in \Delta: \#F=d-1,\ \mbox{$F$ is  contained in exactly one facet of $\Delta$}\rangle$$
is a simplicial $(d-2)$-sphere.

\begin{definition}
Let $M$ be a proper $\cc$-multicomplex.
We call the boundary $\bier_\cc (M)\\
=\partial \B_\cc(M)$ of the $(|\cc|-1)$-dimensional simplicial ball $\B_\cc(M)$ the
\textit{Bier sphere of a multicomplex $M$ with respect to $\cc$}.
\end{definition}

In the rest of this section,
we study some easy combinatorial properties of Bier spheres.
First,
we describe the facets of $\bier_\cc (M)$.
For any monomial $x^\aaa$ and for any pure power of a variable $x_i^j$,
we define
$$x^\aaa \diamond x_i^j= x^\aaa x_i^{-a_i+j} = x_1^{a_1} \cdots x_{i-1}^{a_{i-1}} x_i^j x_{i+1}^{a_{i+1}} \cdots x_n^{a_n}$$
(we assume $x_i^0 \ne 1$ when we consider $x^\aaa \diamond x_i^0$).
The facets of $\bier_\cc (M)$ are given as follows.

\begin{proposition} \label{facets}
Let $M$ be a proper $\cc$-multicomplex. Then
$$\bier_\cc(M)=\big\langle 
F_\cc(x^\aaa)\setminus \{x_i^{(j)} \}: x^\aaa \in M,\ x^\aaa \diamond x_i^j \not \in M,\
a_i<j \leq c_i \big \rangle.
$$
\end{proposition}

\begin{proof}
Let $x^\aaa \in M$.
Note that $x^\aaa \diamond x_i^j \not \in M$ implies $a_i<j$.
$\bier_\cc(M)$ is generated by all codimension $1$ faces $F$ of $\B_\cc(M)$ such that $F$ is contained in  exactly one facet of $\B_\cc(M)$.
It is easy to see that if $F_\cc(x^\aaa) \setminus \{x_i^{(j)}\}$,
where $x_i^{(j)} \in F_\cc(x^\aaa)$,
is contained in
$F_\cc(x^\bb)$ for some $\cc$-monomial $x^\bb$, then $x^\bb$ must be either
$x^\aaa$ or $x^\aaa \diamond x_i^j$.
This implies the desired formula.
\end{proof}

\begin{example}
Let $\cc=(2,2)$ and $M=\{1,x,y,x^2,y^2\}$. Then
\begin{align*}
\B_{(2,2)}(M)=
\left \langle
\begin{array}{l}
x^{(1)}x^{(2)}y^{(1)}y^{(2)},\
x^{(0)}x^{(2)}y^{(1)}y^{(2)},\
x^{(1)}x^{(2)}y^{(0)}y^{(2)},\\
x^{(0)}x^{(1)}y^{(1)}y^{(2)},\
x^{(1)}x^{(2)}y^{(0)}y^{(1)}
\end{array}
\right \rangle,
\end{align*}
where we identify $x^{(i)}x^{(j)}y^{(k)}y^{(\ell)}$ with $\{x^{(i)},x^{(j)},y^{(k)},y^{(\ell)}\}$ for simplicity, and
\begin{align*}
\bier_{(2,2)}(M)=
\left \langle
\begin{array}{l}
x^{(0)}x^{(2)}y^{(1)},\
x^{(0)}x^{(2)}y^{(2)},\
x^{(1)}y^{(0)}y^{(2)},\
x^{(2)}y^{(0)}y^{(2)},\\
x^{(0)}x^{(1)}y^{(1)},\
x^{(0)}x^{(1)}y^{(2)},\
x^{(1)}y^{(0)}y^{(1)},\
x^{(2)}y^{(0)}y^{(1)}
\end{array}
\right \rangle.
\end{align*}
Then $\bier_{(2,2)}(M)$ is the boundary complex of the octahedron, that is,
$\bier_{(2,2)}(M)=\partial \langle \{x^{(0)},y^{(0)}\} \rangle*\partial \langle \{x^{(1)},x^{(2)}\}\rangle*\partial \langle \{y^{(1)},y^{(2)}\}\rangle$.
\end{example}

\begin{example}
\label{codim1}
Here we classify Bier spheres of multicomplexes with one variable.
Let $\cc=(c)$ and $M=\{1,x,x^2,\dots,x^b\}$, where $0 \leq b < c$.
Then
$$
\B_{(c)}(M)=
\left \langle
\{x^{(0)},x^{(1)}, \cdots, x^{(c)}\} \setminus \{x^{(i)}\}: i=0,1,\dots,b
\right \rangle.
$$
and
$$
\ \bier_{(c)}(M)=
\left \langle
\{x^{(0)},x^{(1)}, \cdots, x^{(c)}\} \setminus \{x^{(i)},x^{(j)}\}: i=0,1,\dots,b,\ j=b+1,\dots,c
\right \rangle.
$$
Hence $\bier_{(c)}(M)=\partial \langle \{x^{(0)},x^{(1)}, \cdots, x^{(b)}\}\rangle * \partial \langle \{x^{(b+1)},x^{(b+2)}, \cdots, x^{(c)}\}\rangle$.
\end{example}

If $\cc=(1,1,\dots,1)$ then $\cc$-multicomplexes can be identified with simplicial complexes.
We show that $\bier_{(1,1,\dots,1)}(M)$ are classical Bier spheres.
Let $\Delta$ and $\Gamma$ be  simplicial complexes on the vertex sets $V$ and $W$, respectively.
We say that $\Delta$ is isomorphic to $\Gamma$ if  there is a bijection $\phi: V \to W$
such that $\Gamma=\{ \phi(F):F \in \Delta\}$.

\begin{theorem} \label{bier}
Let $\cc=(1,1,\dots,1)$ and $M$ a proper $\cc$-multicomplex on $X$.
Let $\Delta$ be the simplicial complex defined by $\Delta=\{\{i_1,\dots,i_k\} \subset [n]: x_{i_1} \cdots x_{i_k} \in M\}$.
Then $\bier_\cc(M)$ is combinatorially isomorphic to $\bier(\Delta)$.
\end{theorem}

\begin{proof}
A straightforward computation shows (see \cite[Lemma 5.6.4]{Ma})
$$\bier (\Delta) =
\left\langle X_F \cup Y_{[n]\setminus (F\cup\{j\})}:
F \in \Delta,\ F\cup \{j\} \not \in \Delta \right \rangle.
$$
For $F \subset [n]$, let $x^F = \prod_{i \in F} x_i$.
Observe that $F_\cc(x^F)=\{x_i^{(0)}:i \in F\}\cup \{x_i^{(1)}:i \not\in F\}$.
Proposition \ref{facets} shows
$$
\bier_\cc(M)=
\left\langle
F_\cc(x^F)\setminus\{x_j^{(1)}\}:F \in \Delta,\ F\cup\{j\} \not \in \Delta
\right \rangle,
$$
which proves the desired statement.
\end{proof}

Lemma \ref{ball} and Theorem \ref{bier}
provide a new proof of Bier's theorem.
The known proofs of Bier's theorem uses subdivisions or Bistellar flips.
See \cite{BPSZ,de,Ma}.

Finally, we compute face vectors of Bier spheres.
Let $\Delta$ be a $(d-1)$-dimensional simplicial complex.
Let $f_i=f_i(\Delta)$ be the number of faces of $\Delta$ of cardinality $i$.
The \textit{$f$-vector} of $\Delta$ is the vector $f(\Delta)=(f_0,f_1,\dots,f_d)$ where $f_0=1$,
and the \textit{$h$-vector} $h(\Delta)=(h_0,h_1,\dots,h_d) \in \ZZ^{d+1}$ of $\Delta$ is defined by the relation
$$\sum_{i=0}^d h_i t^{d-i} = \sum_{i=0}^d f_i (t-1)^{d-i}.$$
Also, for $0 \leq i \leq \lfloor \frac d 2 \rfloor$, let $g_i= h_i -h_{i-1}$ where $g_0=1$.
The vector $g(\Delta)=(g_0,g_1,\dots,g_{\lfloor \frac d 2 \rfloor})$ is called the \textit{$g$-vector of $\Delta$}.
For a multicomplex $M$,
let $f_i(M)$ be the number of monomials in $M$ of degree $i$.

It is easy to see that knowing the $f$-vector of $\Delta$ is equivalent to knowing the $h$-vector of $\Delta$.
Also, if $\Delta$ is a simplicial sphere then knowing the $g$-vector of $\Delta$ is equivalent to knowing the $h$-vector of $\Delta$
by the Dehn-Sommerville equations $h_i=h_{d-i}$.
The $g$-vectors of Bier spheres are given by the following formula.

\begin{theorem}
\label{hvector}
Let $M$ be a proper $\cc$-multicomplex.
\begin{itemize}
\item[(i)] $h_i (\B_\cc (M))= f_i(M)$ for all $i$.
\item[(ii)] $g_i(\bier_\cc (M))= f_i(M) - f_{|\cc|-i}(M)$ for $i=0,1,\dots, \lfloor \frac {|\cc|-1} 2 \rfloor$.
\end{itemize}
\end{theorem}

\begin{proof}
(i) follows from the shelling \eqref{1.1} and the following fact:
If $\Delta$ is shellable with a shelling $F_1,\dots,F_s$ then, for $i \geq 1$,
\begin{align}
\label{shell}
h_i(\Delta)=\#\{k: k\geq 2,\ \langle F_1,\dots,F_{k-1}\rangle \cap \langle F_k \rangle
\mbox{ is generated by $i$ facets}\}.
\end{align}
See \cite[III, Proposition 2.3]{St}.

(ii) follows from the well-known fact that
if $B$ is a $(d-1)$-dimensional simplicial ball then $g_i(\partial B) = h_i(B)-h_{d-i}(B)$
for $i=0,1,\dots,\lfloor \frac {d-1} 2 \rfloor$.
See \cite[p.\ 137]{St}.
\end{proof}

\section{Shellability}

In the previous section, we prove that $\B_\cc(M)$ is a simplicial ball by using shellability.
Shellability is an important property in combinatorial topology,
which induces several important topological and enumerative properties
like Lemma \ref{shellball} and formula \eqref{shell}.
Thus if one obtains a construction of simplicial spheres
it is a fundamental question to ask whether they are shellable.
Bj\"orner et al.\ \cite{BPSZ} proved that Bier spheres of simplicial complexes are shellable.
In this section,
we extend this result for Bier spheres of multicomplexes.

Let $\cc=(c_1,\dots,c_n) \in \ZZ_{\geq 0}^n$ and let $M$ be a proper $\cc$-multicomplex.
For any $\cc$-monomial $x^\aaa$ and for any pure power of a variable $x_i^j$ with $0 < j \leq c_i$,
let
$$G(x^\aaa;x_i^j)=F_\cc(x^\aaa) \setminus \{x_i^{(j)}\},$$
and let
$$
\facet (M) =
\big\{ G(x^\aaa;x_i^j) : x^\aaa \in M,\ x^\aaa \diamond x_i^j \not \in M,\ a_i < j \leq c_i \big\}.
$$
By Proposition \ref{facets},
$\facet (M)$ is the set of the facets of $\bier_\cc (M)$.

We introduce the order $\succ$ on $\facet (M)$ as follows:
Let $>_\lex$ be the lexicographic order on $\ZZ^n$.
Thus, for $\aaa, \bb \in \ZZ^n$,
$\aaa >_\lex \bb$ if and only if there exists an $i$ such that $a_i > b_i$ and $a_k = b_k$ for all $k<i$.
Let $G(x^\aaa;x_p^s), G(x^\bb;x_q^t) \in \facet (M)$.
We define $G(x^\aaa;x_p^s) \succ G(x^\bb;x_q^t)$ if one of the following conditions hold
\begin{itemize}
\item[(i)]
$(a_1,a_2,\dots,a_{p-1},s,-a_{p+1},\dots,-a_n) >_\lex (b_1,b_2,\dots,b_{q-1},t,-b_{q+1},\dots,-b_n),$
\item[(ii)] $(a_1,a_2,\dots,a_{p-1},s,-a_{p+1},\dots,-a_n) = (b_1,b_2,\dots,b_{q-1},t,-b_{q+1},\dots,-b_n)$ and $x_p^{a_p+1}>_\lex x_q^{b_q+1}$,
\end{itemize}
where we define $x_p^{a_p+1}>_\lex x_q^{b_q+1}$ if $p<q$ or $p=q$ and $a_p>b_q$.
Clearly, $\succ$ is a total order on $\facet (M)$.

\begin{theorem} \label{shellable}
For any proper $\cc$-multicomplex $M$,
$\bier_\cc (M)$ is shellable.
\end{theorem}

\begin{proof}
We show that the order $\succ$ on $\facet (M)$ induces a shelling of $\bier_\cc(M)$.
Fix $G(x^\aaa;x_p^s) \in \facet (M)$.
Let
$$\Gamma = \big\langle G(x^\bb;x_q^t) \in \facet (M): G(x^\bb;x_q^t) \succ G(x^\aaa;x_p^s)\big\rangle,$$
and let
\begin{align*}
H=&\ \{x_i^{(j)}\in \widetilde X: i<p \mbox{ and } j>a_i\}\\
&\ \bigcup \{x_i^{(j)}\in \widetilde X: i>p \mbox{ and } j<a_i\}\\
&\ \bigcup \left(\{x_p^{(j)}\in \widetilde X: j>s\} \cup \{x_p^{(j)}\in \widetilde X:  x^\aaa \diamond x_p^j \in M \mbox{ and } j> a_p\} \right).
\end{align*}
Since $s > a_p$, we have
$$H \subset G(x^\aaa;x_p^s)=\widetilde X \setminus \{ x_1^{(a_1)}, \cdots,x_n^{(a_n)},x_p^{(s)}\}.$$
To prove the statement,
it is enough to prove that
\begin{align}
\label{two-one}
\Gamma \cap \big\langle G(x^\aaa;x_p^s) \big\rangle = \big\langle G(x^\aaa;x_p^s) \setminus \{ x_i^{(j)}\} : x_i^{(j)} \in H \big\rangle.
\end{align}

Let $u=x^\aaa \diamond x_p^s$.
The inclusion `$\supset$' follows from the following case analysis.\smallskip

[\textit{Case 1.1}]
Suppose $i<p$ and $x^\aaa \diamond x_i^j \not \in M$ with $j > a_i$.
Then $G(x^\aaa;x_i^j)\in \facet(M)$ satisfies
$G(x^\aaa;x_i^j) \succ G(x^\aaa;x_p^s)$ and
$G(x^\aaa;x_i^j) \supset G(x^\aaa;x_p^s) \setminus \{x_i^{(j)}\}$.

[\textit{Case 1.2}]
Suppose $i<p$ and $x^\aaa \diamond x_i^j \in M$ with $j > a_i$.
Then $G(x^\aaa\diamond x_i^j ;x_p^s)\in \facet(M)$ satisfies
$G(x^\aaa\diamond x_i^j ;x_p^s) \succ G(x^\aaa;x_p^s)$ and
$G(x^\aaa\diamond x_i^j ;x_p^s) \supset G(x^\aaa;x_p^s) \setminus \{x_i^{(j)}\}$.

[\textit{Case 2.1}]
Suppose $i>p$ and $u \diamond x_i^j \not \in M$ with $j < a_i$.
Then $G(x^\aaa\diamond x_i^j;x_p^s)\in \facet(M)$ satisfies
$G(x^\aaa\diamond x_i^j;x_p^s) \succ G(x^\aaa;x_p^s)$ and
$G(x^\aaa\diamond x_i^j;x_p^s) \supset G(x^\aaa;x_p^s) \setminus \{x_i^{(j)}\}$.

[\textit{Case 2.2}]
Suppose $i>p$ and $u \diamond x_i^j \in M$ with $j < a_i$.
Then $G(u \diamond x_i^j ;x_i^{a_i})\in \facet(M)$ satisfies
$G(u \diamond x_i^j ;x_i^{a_i}) \succ G(x^\aaa;x_p^s)$ and
$G(u \diamond x_i^j ;x_i^{a_i}) \supset G(x^\aaa;x_p^s) \setminus \{x_i^{(j)}\}$.

[\textit{Case 3.1}]
Suppose $i=p$ and  $j > s$.
Then $G(x^\aaa;x_p^j)\in \facet(M)$ satisfies
$G(x^\aaa;x_p^j) \succ G(x^\aaa;x_p^s)$ and
$G(x^\aaa;x_p^j) \supset G(x^\aaa;x_p^s) \setminus \{x_i^{(j)}\}$.

[\textit{Case 3.2}]
Suppose $i=p$ and  $x^\aaa \diamond x_p^j \in M$ with $j>a_p$.
Then $G(x^\aaa\diamond x_p^j;x_p^s)\in \facet(M)$ satisfies
$G(x^\aaa\diamond x_p^j;x_p^s) \succ G(x^\aaa;x_p^s)$ and
$G(x^\aaa\diamond x_p^j;x_p^s) \supset G(x^\aaa;x_p^s) \setminus \{x_i^{(j)}\}$.\medskip

Next, we prove that the right-hand side contains the left-hand side in \eqref{two-one}.
In the same way as in the proof of Lemma \ref{ball}(ii),
what we must prove is $H \not \in \Gamma$.
Suppose $G(x^\bb;x_q^t) \in \facet(M)$ contains $H$ and satisfies $G(x^\bb;x_q^t) \succeq G(x^\aaa;x_p^s)$.
We prove that $G(x^\bb;x_q^t) = G(x^\aaa;x_p^s)$.

Note that $x^\bb \in M$ and $x^\aaa\diamond x_p^s \not \in M$.
Let
$$\aaa'=(a_1,\dots,a_{p-1},s,-a_{p+1},\dots,-a_n)$$
and
$$\mathbf d = (d_1,\dots,d_n) =(b_1,\dots,b_{q-1},t,-b_{q+1},\dots,-b_n).$$
Since $G(x^\bb;x_q^t) \supset H$,
by the definition of $H$, we have
\begin{align}
\label{2.2}
d_1 \leq a_1, \dots, d_{p-1} \leq a_{p-1}, d_p \leq s \mbox{ and } b_{p+1} \geq a_{p+1}, \dots, b_n \geq a_n.
\end{align}
Observe $\mathbf d \geq_\lex \aaa'$ since $G(x^\bb;x_q^t) \succeq G(x^\aaa;x_p^s)$.
If $d_i <a_i$ for some $1 \leq i \leq p-1$ or $d_p <s$ then $\mathbf d <_\lex \aaa'$.
Hence we have
$$
d_1=a_1,\dots,d_{p-1}=a_{p-1} \mbox{ and } d_p =s.
$$
In particular, we have $q \geq p$ since if $q<p$ then $d_p \leq 0$ but $s>a_p \geq 0$ is positive.
We show $q=p$.
If $q>p$ then $b_p=d_p=s$ and $x^\aaa \diamond x_p^s \not \in M$ divides $x^\bb=x_1^{a_1} \cdots x_{p-1}^{a_{p-1}}x_p^s x_{p+1}^{b_{p+1}} \cdots x_n^{b_n}$ by \eqref{2.2}, which contradicts the fact that $x^\bb \in M$.

Now we know $b_1=a_1,\dots,b_{p-1}=a_{p-1}$, $p=q$
and $t=d_p=s$.
By \eqref{2.2} if $b_i > a_i$ for some $p+1 \leq i \leq n$ then $\mathbf d <_\lex \aaa'$.
Thus we have $b_{p+1}=a_{p+1}, \dots,b_n=a_n$,
and therefore $\aaa'=\mathbf d$.

It remains to prove $b_p=a_p$.
Since $G(x^\bb;x_q^t) \succ G(x^\aaa;x_p^s)$,
by the definition (ii) of the order $\succ$ we have $b_p \geq a_p$.
On the other hand, the definition of $H$ says that, for any $x_p^{(j)} \not \in H$ with $j>a_p$, one has $x^\aaa \diamond x_p^j \not \in M$.
This shows that if $b_p>a_p$ then $x^\bb=x^\aaa \diamond x_p^{b_p} \not \in M$ since $G(x^\bb;x_q^t)$ contains $H$,
which contradicts $x^\bb \in M$.
Hence $b_p=a_p$.
\end{proof}

\begin{example}
Let $\cc=(2,2)$ and $M=\{1,x,y,y^2\}$.
Then $\bier_\cc(M)$ is generated by
\begin{align*}
&\big\{
G(1;x^2),G(x;x^2),G(x;y),G(x;y^2),G(y;x),G(y;x^2),G(y^2;x),G(y^2;x^2)\big\}
\\
&=
\left\{
\begin{array}{ll}
x^{(1)}y^{(1)}y^{(2)},x^{(0)}y^{(1)}y^{(2)},x^{(0)}x^{(2)}y^{(2)},
x^{(0)}x^{(2)}y^{(1)},\\
x^{(2)}y^{(0)}y^{(2)},x^{(1)}y^{(0)}y^{(2)},
x^{(2)}y^{(0)}y^{(1)},x^{(1)}y^{(0)}y^{(1)}
\end{array}
\right\}.
\end{align*}
The proof of Theorem \ref{shellable} shows that
$$
G(x;x^2)\! \succ \!
G(1;x^2)\! \succ \!
G(y;x^2)\! \succ \!
G(y^2;x^2)\! \succ \!
G(x;y^2)\! \succ \!
G(x;y)\! \succ \!
G(y;x)\! \succ \!
G(y^2;x)
$$
is a shelling of $\bier_\cc(M)$.
More precisely, the above shelling is
\begin{align*}
& x^{(0)}y^{(1)} y^{(2)} \succ 
\dot x^{(1)}y^{(1)}y^{(2)} \succ 
x^{(1)}\dot y^{(0)}y^{(2)} \succ 
x^{(1)}\dot y^{(0)} \dot y^{(1)}\! \succ 
x^{(0)}\dot x^{(2)}y^{(1)}
 \succ
x^{(0)}\dot x^{(2)} \dot y^{(2)} \\
&\succ 
\dot x^{(2)} \dot y^{(0)}y^{(2)} \succ 
\dot x^{(2)} \dot y^{(0)} \dot y^{(1)},
\end{align*}
where variables with $\cdot$ correspond to variables in $H$.
\end{example}

\section{Bier spheres and polarization}

The construction of $\bier_\cc (M)$ is inspired by the study of polarizations
of monomial ideals in commutative algebra.
In this section, we study connections between Bier spheres and polarization.

We recall some basics on commutative algebra.
Let $S=K[x_1,\dots,x_n]$ be the polynomial ring over a field $K$ with $\deg x_i=1$ for any $i$.
For a simplicial complex $\Delta$ on $X=\{x_1,\dots,x_n\}$,
the ideal
$$I_\Delta=(x_{i_1}\cdots x_{i_k}: \{x_{i_1},\dots,x_{i_k}\} \subset X,\ \{x_{i_1},\dots,x_{i_k}\} \not \in \Delta)$$
is called the \textit{Stanley-Reisner ideal of $\Delta$}.
The ring
$$K[\Delta]=S/I_\Delta$$
is called the \textit{Stanley-Reisner ring of $\Delta$}.
Note that the correspondence $\Delta \leftrightarrow I_\Delta$ gives a one-to-one correspondence between
simplicial complexes on $X$ and squarefree monomial ideals in $S$.

Let $I \subset S$ be a homogeneous ideal and $R=S/I$.
The \textit{(Krull) dimension} of $R$ is the maximal number of homogeneous elements of $R$ which are algebraically independent over $K$.
The \textit{Hilbert series of $R$} is the formal power series $H_R(t)=\sum_{k \geq 0} (\dim_K R_k) t^k$,
where $R_k$ is the graded component of $R$ of degree $k$.
It is known that $H_R(t)$ is a rational function of the form
$H_R(t)=(h_0+h_1t+\cdots+h_s t^s)/ {(1-t)^d},$
where $h_s \ne 0$ and where $d$ is the Krull dimension of $R$
\cite[Corollary 4.1.8]{BH}.
We write $h_i(R)=h_i$ for all $i \geq 0$,
where $h_i=0$ if $i >s$.
The vector $(h_0,h_1,\dots,h_s)$ is called the \textit{$h$-vector of $R$}.
Note that if $R=K[\Delta]$ then $s \leq d$, $h_i(\Delta)=h_i(R)$ for all $i \leq d$
and the Krull dimension of $R$ is equal to $\dim \Delta+1$.

Fix $\cc=(c_1,c_2,\dots,c_n) \in \ZZ_{\geq 0}^n$.
Let $\overline \cc=(c_1+1,c_2+1,\dots,c_n+1)$.
For a monomial ideal $I \subset S$,
we write $G(I)$ for the unique minimal set of monomial generators of $I$.
A monomial ideal $I \subset S$ is called a \textit{$\overline \cc$-ideal}
if $I$ is generated by $\overline \cc$-monomials.

\begin{definition}
Let $I \subset S$ be a $\overline \cc$-ideal
and let $\widetilde S=K[x_{i,j}:1 \leq i \leq n, 0 \leq j \leq c_i]$.
The {\em polarization of a $\overline \cc$-monomial $x^\aaa=x_1^{a_1}x_2^{a_2}\cdots x_n^{a_n}$ (with respect to $\overline \cc$)} is the squarefree monomial
$$\pol(x^\aaa)=
\prod_{a_i \ne 0} (x_{i,0} x_{i,1}\cdots x_{i,a_i-1}) \in \widetilde S.$$
The \textit{polarization of $I$ (with respect to $\overline \cc$)} is the squarefree monomial ideal
$$\pol(I)=\big(\pol(x^\aaa): x^\aaa \in G(I)\big) \subset \widetilde S.$$
\end{definition}

The ideal $\pol(I)$ has following properties (see e.g., \cite[pp.\ 59--60]{MS}):
\begin{itemize}
\item $S/I$ and $\widetilde S/\pol(I)$ have the same graded Betti numbers,
in particular, have the same $h$-vector.
\item $\dim \widetilde S/\pol(I)= \dim S/I+ |\cc|$.
\end{itemize}
In the rest of this section,
we identify $x_{i,j}$ and $x_i^{(j)}$,
and regard $\B_\cc(M)$ and $\bier_\cc(M)$ as simplicial complexes on $V=\{x_{i,j}:1\leq i\leq n, 0 \leq j \leq c_i\}$.
For a $\cc$-multicomplex $M$,
let $I(M) \subset S$ be the ideal generated by all monomials which are not in $M$.
If $M$ is a $\cc$-multicomplex then $I(M)$ is a $\overline \cc$-ideal since $I(M)$ contains $x_1^{c_1+1},\dots,x_n^{c_n+1}$.

\begin{lemma} \label{jahan}
For any $\cc$-multicomplex $M$,
$I_{\B_\cc(M)}=\pol(I(M)).$
\end{lemma}

\begin{proof}
Let $\Delta$ be the simplicial complex such that $I_\Delta=\pol(I(M))$.
We claim $\Delta = \B_\cc(M)$.
Since taking polarization does not change $h$-vectors,
by Lemma \ref{hvector}(i),
$$h_i(\Delta)=h_i(\widetilde S/I_\Delta)=h_i(S/I(M))=
\#\{x^\aaa \in M : \deg x^\aaa=i\}=h_i(\B_\cc(M)).$$
Since $\dim \Delta=\dim(\widetilde S/\pol(I(M)))-1=|\cc|-1$,
the above equations show that $\Delta$ and $\B_\cc(M)$ have the same $f$-vector.

Let $F=F_\cc(x^\aaa)=\widetilde X\setminus \{ x_1^{(a_1)},\dots,x_n^{(a_n)}\}$,
where $x^\aaa \in M$,
be a facet of $\B_\cc(M)$.
Since $\Delta$ and $\B_\cc(M)$ have the same $f$-vector,
to prove $\Delta=\B_\cc(M)$,
it is enough to prove that $\prod_{x_i^{(j)} \in F} x_{i,j}$ is not contained in $I_\Delta=\pol(I(M))$.
Suppose contrary that $\prod_{x_i^{(j)} \in F} x_{i,j} \in \pol(I(M))$.
Then there is a monomial $x^\bb \in I(M)$
such that $\pol(x^\bb)$ divides $\prod_{x_i^{(j)} \in F} x_{i,j}$.
By the definition of polarization, we have $b_1 \leq a_1,\dots,b_n \leq a_n$.
Then $x^\bb$ divides $x^\aaa$, which contradicts $x^\aaa \not \in I(M)$.
\end{proof}

\begin{remark}
Jahan \cite[Lemma~3.7]{So} computed the facets of the simplicial complex whose Stanley-Reisner ideal is the polarization of a monomial ideal.
This will give an alternative proof of Lemma \ref{jahan}.
See also \cite[Theorem 3.1]{HP}.
\end{remark}

Lemma \ref{jahan} allows us to study the structure of $\bier_\cc(M)$ from an algebraic viewpoint.
We discuss some algebraic aspects of Bier spheres later in Section 5.

Next we study generators of $I_{\bier_\cc(M)}$.
Let $\Delta$ be a simplicial complex on $[n]$.
Recall that the Alexander dual of $\Delta$
is the simplicial complex $\Delta^\vee=\{F \subset [n]: [n] \setminus F \not \in \Delta\}$.
Let $T=K[x_1,\dots,x_n,y_1,\dots,y_n]$,
$X_\Delta =(x^F: F \not \in \Delta) \subset T$
and $Y_{\Delta^\vee} =(y^F: F \not \in \Delta^\vee) \subset T$,
where $x^F=\prod_{i \in F} x_i$ and $y^F=\prod_{i \in F} y_i$.
Then it is not hard to see that
$$I_{\bier(\Delta)}
=X_\Delta  + Y_{\Delta^\vee} +(x_1y_1,x_2y_2,\dots,x_ny_n).$$
Indeed, 
one has $x^F y^G \not \in I_{\bier(\Delta)}$
if and only if $F \cap G = \emptyset$,
$x^F \not \in X_\Delta$ and $y^G \not \in Y_{\Delta^\vee}$,
which implies the above equation.
We show that a similar formula holds for
Bier spheres of multicomplexes.
We first recall Alexander duality of monomial ideals
introduced by Miller \cite[Definition 1.5]{Mi2}.

\begin{definition}
\label{alex}
Let $M$ be a $\cc$-multicomplex.
The \textit{Alexander dual of $M$ with respect to $\cc$} is the multicomplex defined by
$$
M^{\vee}=\{ x_1^{c_1-a_1}\cdots x_n^{c_n-a_n}: x^\aaa \mbox{ is a $\cc $-monomial such that } x^\aaa \not \in M\}.
$$
Note that $M^\vee$ depends not only on $M$ but also on $\cc$.

We also write an ideal-theoretic definition of Alexander duality.
For a $\cc$-ideal $I \subset S$,
we call the ideal
$$I^\vee=(x^{c_1-a_1}\cdots x_n^{c_n-a_n}: x^\aaa \mbox{ is a $\cc $-monomial such that } x^\aaa \not \in I) \subset S$$
the \textit{Alexander dual of $I$ with respect to $\cc$.}
\end{definition}

Let $I_\cc(M) \subset S$ be the monomial ideal generated by all $\cc$-monomials which are not in $M$.
Then the above two definitions are related by $I_\cc(M)^\vee = I_\cc(M^\vee)$.

For any $\overline \cc$-monomial $x^\aaa \in S$,
let
$$\pol^{*}(x^\aaa)= \prod_{a_i \ne 0} (x_{i,c_i}x_{i,c_i-1}\cdots x_{i,c_i-a_i+1}).$$
Similarly, for a $\overline \cc$-ideal $I \subset S$, let $\pol^* (I) \subset \widetilde S$ be the ideal generated by $\{\pol^* (x^\aaa): x^\aaa \in G(I)\}$.
(Thus $\pol^*$(-) is the polarization by using an opposite ordering of the variables $x_{i,c_i},x_{i,c_i-1},\dots,x_{i,1}$.)

\begin{lemma}
\label{complement}
Let $M$ be a proper $\cc$-multicomplex and let
$$\mathcal B^*_\cc(M)= \langle F_\cc(x^\aaa): x^\aaa\mbox{ is a $\cc$-monomial such that } x^\aaa \not \in M\rangle$$
be the complementary ball of $\B_\cc(M)$ in $\Lambda_\cc$.
Then
$I_{\mathcal B^*_\cc(M)}= \pol^*(I(M^\vee))$.
\end{lemma}

\begin{proof}
Let $\pi$ be the permutation on the vertex set $V$ defined by $\pi(x_{i,j})=x_{i,c_i-j}$.
Then $\pi(F_\cc(x^\aaa))=F_\cc(x_1^{c_1-a_1}\cdots x_n^{c_n-a_n})$ and
$\pi(\B_\cc(M^\vee))=\B_\cc^*(M)$.
Hence
$$I_{\mathcal B^*_\cc(M)}= \pi (I_{\B_\cc(M^\vee)})= \pi \circ \pol(I(M^\vee))=\pol^*(I(M^\vee)),$$
as desired.
\end{proof}

\begin{theorem}
\label{generator}
Let $M$ be a proper $\cc$-multicomplex.
$$ I_{\bier_\cc (M)}=
\pol(I_\cc(M))+ \pol^*(I_\cc(M^{\vee})) +\pol(x_1^{c_1+1},\dots,x_n^{c_n+1}).$$
\end{theorem}

\begin{proof}
Since $\B_\cc^*(M)$ is the complementary ball of $\B_\cc(M)$ in $\Lambda_\cc$,
we have $\B_\cc(M) \cap \B_\cc^*(M)=\partial \B_\cc(M)=\bier_\cc(M)$.
Since $I(M)=I_\cc(M)+(x_1^{c_1+1},\dots,x_n^{c_n+1})$,
by Lemma \ref{complement} we have
$$I_{\bier_\cc(M)} = I_{\B_\cc(M)}+I_{\mathcal B^*_\cc(M)}=
\pol(I_\cc(M))+ \pol^{*}(I_\cc(M^{\vee})) +\pol(x_1^{c_1+1},\dots,x_n^{c_n+1}),$$
as desired.
\end{proof}

\begin{example}
Let $\cc=(2,2,2)$ and $M=\{1,x,y,z,xz,yz,z^2,yz^2\}$.
Then $I_\cc(M)=(x^2,y^2,xy,xz^2)$
and
$I_\cc(M^\vee)=(x y^2 z, x^2 y)$.
Thus
$$
I_{\bier_\cc(M)}=
(x_0x_1,y_0y_1,x_0y_0,x_0z_0z_1)
+(x_2y_2y_1z_2,x_2x_1y_2)+(x_0x_1x_2,y_0y_1y_2,z_0z_1z_2).$$
\end{example}

\begin{example} \label{3..4}
Theorem \ref{generator} may not give minimal generators.
For example, if $\Delta=\langle \{1,2\},\{3\}\rangle$ then
$$I_{\bier(\Delta)}=(x_1x_3,x_2x_3)+(y_3,y_1y_2)+(x_1y_1,x_2y_2,x_3y_3).$$
However, the set of generators $\{x_1x_3,x_2x_3,y_3,y_1y_2,x_1y_1,x_2y_2,x_3y_3\}$
is not minimal.
\end{example}

Finally we note the following result which immediately follows from the fact $(M^\vee)^\vee =M$.

\begin{corollary}
\label{dual}
Let $M$ be a proper $\cc$-multicomplex.
Then $\bier_\cc(M)$ and $\bier_\cc(M^\vee)$ are combinatorially isomorphic.
The isomorphism is given by the permutation of the vertices $x_{i,j} \to x_{i,c_i-j}$ for all $i,j$.
\end{corollary}

\section{Edge decomposability}

In this section, we show that Bier spheres are edge decomposable.
We first recall the definition of edge decomposability.

Let $\Delta$ be a simplicial complex on $V$.
The \textit{link of $\Delta$ with respect to a face $F \in \Delta$}
is the simplicial complex
$$
\lk_\Delta(F) = \{ G \subset V \setminus F : G \cup F \in \Delta\}.
$$
The \textit{contraction $\mathcal C_\Delta (i,j)$ of $\Delta$ with respect to an edge $\{i,j\} \in \Delta$}
is the simplicial complex which is obtained
from $\Delta$ by identifying the vertices $i$ and $j$,
in other words,
$$\mathcal C_\Delta(i,j) =
\{ F \in \Delta: i \not \in F\} \cup \{ (F \setminus \{i\}) \cup \{j\}:
i \in F \in \Delta\}.
$$
(We consider that $\mathcal C_\Delta(i,j)$ is a simplicial complex on $V \setminus \{i\}$.)
We say that $\Delta$ satisfies the \textit{Link condition
with respect to $\{i,j\} \subset V$} if 
\begin{align*}
\lk_\Delta(\{i\}) \cap \lk_\Delta(\{j\})  = \lk_\Delta(\{i,j\}).
\end{align*}

\begin{definition}
\label{edgedecomp}
The boundary complex of a simplex and $\{\emptyset\}$ are \textit{edge decomposable},
and, recursively, a pure simplicial complex $\Delta$ is said to be \textit{edge decomposable}
if there exists an edge $\{i,j\} \in \Delta$ such that $\Delta$ satisfies the Link condition w.r.t.\
$\{i,j\}$ and both $\lk_\Delta(\{i,j\})$ and $\mathcal C_\Delta(i,j)$ are edge decomposable.
\end{definition}

Edge decomposability was introduced by Nevo \cite{Ne} in the study of the $g$-conjecture for spheres,
which states that the $g$-vector of a simplicial sphere is the $f$-vector of a multicomplex.
He proved that the $g$-vector of an edge decomposable sphere is non-negative in \cite{Ne}.
Later, it was proved in \cite{BN,Mu2} that the $g$-vector of an edge decomposable complex
is the $f$-vector of a multicomplex.

Unfortunately, not all spheres are edge decomposable.
For example, the boundary complex of
Lockeberg's non-vertex decomposable $4$-polytope (see \cite{Ha}) does not satisfy the Link condition w.r.t.\ any edge,
and therefore is not edge decomposable.
On the other hand, it was proved in \cite[Proposition 5.4]{Mu2} that Kalai's squeezed spheres are edge decomposable.
This shows that there are many edge decomposable spheres which are not realizable as polytopes.

In the rest of this section,
we prove that Bier spheres are edge decomposable.

\begin{lemma} \label{join}
If $\Delta$ and $\Gamma$ are edge decomposable then $\Delta*\Gamma$ is edge decomposable.
\end{lemma}

\begin{proof}
By the definition of edge decomposability,
we may assume that $\Delta$ and $\Gamma$ are the boundaries of simplexes of dimension at least $1$.
Suppose $\Delta*\Gamma=\partial F * \partial G$, where $F$ and $G$ are simplexes.
Then, for any pair $\{u,v\}$ of vertices, where $u \in F$ and $v \in G$,
it is easy to see that $\Delta*\Gamma$ satisfies the Link condition w.r.t.\ $\{u,v\}$,
$\lk_{\Delta*\Gamma}(\{u,v\})= \partial (F\setminus \{u\}) * \partial (G\setminus \{v\})$ and
$\mathcal C_{\Delta*\Gamma}(u,v)=\partial (F \cup G \setminus \{u\})$.
By using this fact,
the statement follows inductively.
\end{proof}

The above lemma and Example \ref{codim1} show that if $n=1$ then
$\bier_\cc (M)$ is edge decomposable.
We study the case when $n \geq 2$.

\begin{lemma}
\label{reduction}
Let $n \geq 2$ and let $M$ be a proper $\cc$-multicomplex.
If either $\{x_n^{(0)}\}$ or $\{x_n^{(c_n)}\}$ is not in $\bier_\cc(M)$ then there exist $\cc' \in \ZZ_{\geq 0}^{n-1}$
and a $\cc'$-multicomplex $M'$ on $X \setminus \{x_n\}$ such that $\bier_{\cc'}(M')$ is combinatorially isomorphic to $\bier_\cc(M)$.
\end{lemma}

\begin{proof}
By Corollary \ref{dual}, it is enough to consider the case when $\{x_n^{(0)}\} \not \in \bier_\cc(M)$.
Since $x_n^{(0)}  \in G(I_{\bier_\cc(M)})$, by Theorem \ref{generator}, $x_n \in I_\cc(M)$,
and therefore $M$ contains no monomials which are divisible by $x_n$.
Then $M$ is a $(c_1,c_2,\dots,c_{n-1})$-multicomplex on $\{x_1,\dots,x_{n-1}\}$
and
\begin{align*}
\B_\cc(M)=&\ \B_{(c_1,\dots,c_{n-1},0)}(M)* \langle \{x_n^{(1)},\dots,x_n^{(c_n)}\}\rangle \\
\cong&\  \B_{(c_1,\dots,c_{n-1},0)}(M)*\langle \{x_1^{(c_1+1)},\dots,x_1^{(c_1+c_n)}\}\rangle\\
=&\ \B_{(c_1+c_n,\dots,c_{n-1},0)}(M)\\
=&\ \B_{(c_1+c_n,\dots,c_{n-1})}(M),
\end{align*}
as desired.
\end{proof}

\begin{lemma}
\label{linkcondition}
Let $n \geq 2$ and let $M$ be a proper $\cc$-multicomplex.
Suppose that both  $\{x_1^{(c_1)}\}$ and $\{x_n^{(0)}\}$ are in $\bier_\cc(M)$.
Then $\{x_1^{(c_1)},x_n^{(0)}\} \in \bier_\cc(M)$ and
$\bier_\cc (M)$ satisfies the Link condition with respect to $\{x_1^{(c_1)},x_n^{(0)}\}$.
\end{lemma}

\begin{proof}
Note that $c_1>0$ and $c_n>0$ since $\{x_1^{(c_1)}\}$ and $\{x_n^{(0)}\}$ are in $\bier_\cc(M)$.
It is known that a simplicial complex $\Delta$ on $X$ satisfies the Link condition w.r.t.\ $\{x_i,x_j \} \in \Delta$
if and only if $G(I_\Delta)$ has no monomials which are divisible by $x_i x_j$ \cite[Lemma 2.1]{Mu2}.
By Theorem \ref{generator}, $G(I_{\bier_\cc(M)})$ cannot have monomials which are divisible by $x_1^{(c_1)}x_n^{(0)}$
since $G(\pol(I_\cc(M)))$ contains no monomials which are divisible by $x_1^{(c_1)}$
and $G(\pol^*(I_\cc(M^\vee)))$ contains no monomials which are divisible by $x_n^{(0)}$.
In particular, since $x_1^{(c_1)}, x_n^{(0)}$ and $x_1^{(c_1)}x_n^{(0)}$ are not in $G(I_{\bier_\cc(M)})$, we have $\{x_1^{(c_1)},x_n^{(0)}\} \in \bier_\cc(M)$.
\end{proof}

\begin{lemma}
\label{keylemma}
With the same notation as in Lemma \ref{linkcondition},
\begin{itemize}
\item[(i)]
the link of $\bier_\cc(M)$ with respect to $\{x_1^{(c_1)},x_n^{(0)}\}$
is combinatorially isomorphic to some Bier sphere $\bier_{\cc'} (M')$;
\item[(ii)]
the contraction of $\bier_\cc(M)$ with respect to $\{x_1^{(c_1)},x_n^{(0)}\}$
is combinatorially isomorphic to some Bier sphere $\bier_{\cc'} (M')$ with $\cc' \in \ZZ_{\geq 0}^{n-1}$.
\end{itemize}
\end{lemma}

\begin{proof}
(i)
Let $\cc'=(c'_1,c_2',\dots,c_n')=(c_1-1,c_2,\dots,c_n)$
and $M' = \{ x^\aaa \in M: x^\aaa \mbox{ is a $\cc'$-monomial}\}$.
If $M'$ is $\cc'$-full then by Proposition \ref{facets} all facets of $\bier_\cc(M)$ do not contain $x_1^{(c_1)}$,
which contradicts the assumption $\{x_1^{(c_1)} \} \in \bier_\cc(M)$.
Thus $M$ is not $\cc'$-full.
We claim that
\begin{align}
\label{4.1}
\bier_{\cc'}(M')= \lk_{\bier_\cc(M)}(\{x_1^{(c_1)}\}).
\end{align}
Since both complexes are simplicial spheres having the same dimension,
it is enough to prove that $\bier_{\cc'}(M')\subset \lk_{\bier_\cc(M)}(\{x_1^{(c_1)}\}).$

Let $F$ be a facet of $\bier_{\cc'}(M')$.
Then by Proposition \ref{facets}
there exist $x^\aaa \in M'$ and $x_i^j$ with $j \leq c_i'$ such that
$x^\aaa \diamond x_i^j \not \in M'$ and $F=F_{\cc'}(x^\aaa) \setminus \{x_i^{(j)}\}$.
By the definition of $M'$,
$x^\aaa \in M$ and $x^\aaa \diamond x_i^j  \not \in M$.
Hence by Proposition \ref{facets}
$$F_\cc(x^\aaa) \setminus \{x_i^{(j)}\} = F \cup \{x_1^{(c_1)}\}$$
is a facet of $\bier_\cc(M)$, and therefore $F \in \lk_{\bier_\cc(M)}(\{x_1^{(c_1)}\})$.

Let $M''=(M')^{\vee}$ be the Alexander dual of $M'$ with respect to $\cc'$.
Then, by Corollary \ref{dual},
$\lk_{\bier_{\cc'}(M'')}(\{x_n^{(c_n)}\})$
is combinatorially isomorphic to $\lk_{\bier_{\cc'}(M')}(\{x_n^{(0)}\})=\lk_{\bier_\cc(M)}(\{x_1^{(c_1)},x_n^{(0)}\})$.
On the other hand,
\eqref{4.1} says that $\lk_{\bier_{\cc'}(M'')}(\{x_n^{(c_n)}\})$
is equal to the Bier sphere 
$$\bier_{(c_1-1,c_2,\dots,c_{n-1},c_n-1)}\big(\big\{x^\aaa\in M'': x^\aaa \mbox{ is a $(c_1-1,c_2,\dots,c_{n-1},c_n-1)$-monomial}\big\}\big),$$
as desired.
\medskip

(ii)
Let $\cc'=(c'_1,\dots,c'_{n-1})=(c_1+c_n,c_2,\dots,c_{n-1})$.
For any monomial $x_1^sx_n^tu \in M$, where $u$ is a monomial on $x_2,\dots,x_{n-1}$,
let
\begin{align*}
\sigma(x_1^sx_n^tu)=
\begin{cases}
x_1^{s+t} u, & \mbox{ if $s=c_1$ or $t=0,$}\\
0, & \mbox{ otherwise},
\end{cases}
\end{align*}
and
$$M'=\{\sigma(x^\aaa):x^\aaa \in M\mbox{ and } \sigma(x^\aaa) \ne 0\}.$$
We will show that $\bier_{\cc'}(M')$ is combinatorially isomorphic to $\mathcal C_{\bier_\cc (M)}(x_1^{(c_1)},x_n^{(0)})$.

We first prove that $M'$ is indeed a proper $\cc'$-multicomplex.
Since $x^\cc \not \in M$, $x^{\cc'} \not \in M'$.
It is enough to prove that $M'$ is a multicomplex.
If $x_1^\ell u \in M'$ then one has $x_1^\ell u \in M$ if $\ell \leq c_1$
and $x_1^{c_1} x_n^{\ell-c_1} u \in M$ if $\ell >c_1.$
In the former case
$x_1^\ell u/x_i=\sigma(x_1^\ell u/x_i) \in M'$ for any $x_i$ which divides $x_1^\ell u$.
In the latter case
$x_1^\ell u/x_i=\sigma(x_1^{c_1}x_n^{\ell-c_1} u/x_i) \in M'$ for any $x_i$ with $i \geq 2$ which divides $x_1^\ell u$
and
$x_1^\ell u/x_1=\sigma(x_1^{c_1}x_n^{\ell-c_1} u/x_n) \in M'$.
Hence $M'$ is a multicomplex.

Let $\rho$ be the map from the set of subsets of $\{x_i^{(j)}: 1 \leq i \leq n,\ 0 \leq j \leq c_i\}$
to the set of subsets of $\{x_i^{(j)}: 1 \leq i \leq n-1,\ 0 \leq j \leq c'_i\}$
induced by
\begin{align*}
\rho(x_i^{(j)})=
\begin{cases}
x_i^{(j)}, & \mbox{ if $i \ne n$,}\\
x_1^{(c_1+j)}, & \mbox{ if $i =n$.} 
\end{cases}
\end{align*}
Then, since $\rho(x_1^{(c_1)})=\rho(x_n^{(0)})$,
$\rho(\bier_\cc(M))$ is combinatorially isomorphic to
the contraction $\mathcal C_{\bier_\cc(M)}(x_1^{(c_1)},x_n^{(0)})$.
We claim that
$$\rho(\bier_\cc(M))
=\bier_{\cc'}(M').$$
It follows from \cite[Theorem 1.4]{Ne} and  Lemma \ref{linkcondition} that
$\rho(\bier_\cc(M))$ is a simplicial sphere.
Thus it is enough to prove $\rho(\bier_\cc(M))
\supset \bier_{\cc'}(M')$.

Let $x_1^\ell u \in M'$, where $u$ is a monomial on $x_2,\dots,x_{n-1}$,
and let $x_i^j$ with $j \leq c'_i$ be such that $x_1^\ell u \diamond x_i^j \not \in M'$.
By Proposition \ref{facets},
what we must prove is $F_{\cc'}(x_1^\ell u )\setminus \{x_i^{(j)}\} \in \rho(\bier_\cc(M))$.
This follows from the following case analysis.
\medskip

[\textit{Case 1.1}]
Suppose $i \ne 1$ and $\ell \leq c_1$.
Then $x_1^\ell u \in M$ and $x_1^\ell  u \diamond x_i^j \not \in M$.
$F_{\cc'}(x_1^\ell u )\setminus \{x_i^{(j)}\} = \rho(F_\cc(x_1^\ell u)\setminus \{x_i^{(j)}\}) \in \rho(\bier_\cc(M)).$

[\textit{Case 1.2}]
Suppose $i \ne 1$ and $\ell > c_1$.
Then $x_1^{c_1}x_n^{\ell-c_1} u \in M$ and $x_1^{c_1}x_n^{\ell-c_1}  u \diamond x_i^j \not \in M$.
$F_{\cc'}(x_1^\ell u )\setminus \{x_i^{(j)}\} = \rho(F_\cc(x_1^{c_1}x_n^{\ell-c_1} u)\setminus \{x_i^{(j)}\}) \in \rho(\bier_\cc(M)).$

[\textit{Case 2.1}]
Suppose $i=1$ and $\ell <j  \leq c_1$.
Then $x_1^\ell u \in M$ and $x_1^\ell u \diamond x_1^j \not \in M$.
$F_{\cc'}(x_1^\ell u )\setminus \{x_1^{(j)}\} = \rho(F_\cc(x_1^\ell u)\setminus \{x_1^{(j)}\}) \in \rho(\bier_\cc(M)).$

[\textit{Case 2.2}]
Suppose $i=1$ and $\ell \leq c_1 <j$.
Then $x_1^\ell u \in M$ and $x_1^{c_1} x_n^{j-c_1} u \not \in M$.
If $x_1^\ell x_n^{j-c_1} u \not \in M$ then
$F_{\cc'}(x_1^\ell u )\setminus \{x_1^{(j)}\} = \rho(F_\cc(x_1^\ell u)\setminus \{x_n^{(j-c_1)}\}) \in \rho(\bier_\cc(M)).$
If $x_1^\ell x_n^{j-c_1} u \in M$ then
$F_{\cc'}(x_1^\ell u )\setminus \{x_1^{(j)}\} = \rho(F_\cc(x_1^\ell x_n^{j-c_1} u)\setminus \{x_1^{(c_1)}\}) \in \rho(\bier_\cc(M)).$

[\textit{Case 2.3}]
Suppose $i=1$ and $ c_1 <\ell <j$.
Then $x_1^{c_1}x_n^{\ell-c_1} u \in M$ and $x_1^{c_1} x_n^{j-c_1} u \not \in M$.
$F_{\cc'}(x_1^\ell u )\setminus \{x_i^{(j)}\} = \rho(F_\cc(x_1^{c_1}x_n^{\ell-c_1} u)\setminus \{x_n^{(j-c_1)}\}) \in \rho(\bier_\cc(M)).$
\end{proof}

\begin{theorem}
\label{edgedecomposable}
For any proper $\cc$-multicomplex $M$,
$\bier_\cc(M)$ is edge decomposable.
\end{theorem}

\begin{proof}
We use induction on dimension and $n$.
If $n=1$ then the statement follows from Lemma \ref{join}.
Also, any simplicial $(d-1)$-sphere is edge decomposable if $d \leq 2$.
Suppose $n>1$.
By Lemma \ref{reduction} we may assume that $\{x_1^{(c_1)}\}$ and $\{x_n^{(0)}\}$ are in $\bier_\cc(M)$.
Then the statement follows from Lemmas \ref{linkcondition} and \ref{keylemma} together with the induction hypothesis.
\end{proof}

By \cite[Corollary 3.5]{Mu2},
Theorem \ref{edgedecomposable} implies the following fact.

\begin{corollary} \label{g}
The $g$-vector of $\bier_\cc(M)$ is the $f$-vector of a multicomplex.
\end{corollary}

We proved the above corollary by using edge decomposability.
On the other hand,
since we obtain an explicit formula of the $g$-vector of $\bier_\cc(M)$ in Theorem \ref{hvector}(ii),
it would be desirable to find a simple combinatorial proof of the above corollary
by using that formula.
For Bier spheres of simplicial complexes, Bj\"orner et al.\ \cite[Corollary 5.5]{BPSZ} gave such a proof,
and their proof can be extended to the case when $c_1=\dots=c_n$ by using Clements-Lindstr\"om theorem \cite{CL}.
However, their method seems not to be applicable when $c_1,c_2,\dots,c_n$ are not equal.

Recall that 
any Bier sphere $\bier_\cc(M)$ is the boundary of a simplicial ball $\B_\cc(M)$ which is a subcomplex of the simplicial sphere $\Lambda_\cc
=\partial \widetilde X_1 * \cdots * \partial \widetilde X_n$ with the same dimension.
Corollary \ref{g} can be extended as follows.

\begin{corollary}
If $B \subset \Lambda_\cc$ is a simplicial ball with the same dimension as $\Lambda_\cc$,
then the $g$-vector of $\partial B$ is the $f$-vector of a multicomplex.
\end{corollary}

\begin{proof}
Any subcomplex of $\Lambda_\cc$ with the same dimension as $\Lambda_\cc$ is a balanced complex of type $(c_1-1,\dots,c_n-1)$ \cite{St2}. 
Since $B$ is Cohen-Macaulay,
it follows from \cite[Theorem 4.4]{St2} that the $h$-vector of $B$ is equal to the $h$-vector of a $\cc$-multicomplex $M$.
Then $h(B)=h(\B_\cc(M))$ and $g(\partial B)=g(\bier_\cc(M))$.
\end{proof}


\section{An algebraic study of Bier spheres}


In this section, we study algebraic aspects of Bier spheres.
Fix $\cc=(c_1,\dots,c_n) \in \ZZ_{\geq 1}^n$.
Let $S=K[x_1,\dots,x_n]$ and $\widetilde S=K[x_{i,j}:1\leq i \leq n,\ 0 \leq j \leq c_i]$.

\subsection*{An algebraic proof of Bier's theorem}
We first introduce some basic notations on commutative algebra.
Let $I \subset S$ be a homogeneous ideal and $R=S/I$.
A sequence $\theta_1,\dots,\theta_r \in R$ is said to be an $R$-sequence if $(\theta_1,\dots,\theta_r)R \ne R$ and
$\theta_i$ is not a zero divisor of $R/(\theta_1,\dots,\theta_{i-1})$ for all $i$.
The ring $R$ is said to be \textit{Cohen-Macaulay} if there is an $R$-sequence $\theta_1,\dots,\theta_d$ of length $d$,
where $d$ is the Krull dimension of $R$.
If $R$ is Cohen-Macaulay then the number $\dim_K \{m \in R/(\theta_1,\dots,\theta_d): (x_1,\dots,x_n)m=0\}$ is independent of a choice of an $R$-sequence
$\theta_1,\dots,\theta_d$ and this number is called the \textit{type} of $R$ \cite[I ,Theorem 12.4]{St}.
We say that $R$ is \textit{Gorenstein} if it is a Cohen-Macaulay ring of type $1$.

The Gorenstein property is important in commutative algebra since it implies many nice symmetry.
The Gorenstein property also appears in several combinatorial situations in which symmetry appears.
See e.g., \cite{St3} and \cite[II \S 5]{St}.
For simplicity, we say that an ideal $I \subset S$ is Gorenstein if $S/I$ is Gorenstein.

When we consider the correspondence $\Delta \leftrightarrow K[\Delta]$,
the property that $\Delta$ is a sphere is close to the Gorenstein property of $K[\Delta]$.
Indeed, $K[\Delta]$ is Gorenstein if and only if $\Delta$ is the join of a simplex and a homology sphere
\cite[II Theorem 5.1]{St}.
In particular,
Bier's theorem shows that the ideal in $K[x_1,\dots,x_n,y_1,\dots,y_n]$ generated by
$$\{x^F : F \not \in \Delta\} \cup \{y^F: F \not \in \Delta^\vee\} \cup \{x_1y_1,\dots,x_ny_n\}$$
is Gorenstein.
Although Bier's proof is simple (see \cite[pp.\ 112--116]{Ma}),
it is natural to ask if there is an algebraic way to prove that the above ideal is Gorenstein.
Linkage theory gives such a proof.

We will not explain details on linkage theory,
since we only need the following fact:
Let $Q \subset S$ be a Gorenstein ideal,
$L \subset S$ a Cohen-Macaulay ideal which contains $Q$
and $L'=Q:L=\{f \in S :f L \subset Q\}$.
Suppose that $Q$ and $L$ are radical ideals and $\dim S/L=\dim S/Q$.
Then $S/(L+L')$ is Gorenstein.
(This fact is an immediate consequence of \cite[Theorem 4.2.1 and Proposition 5.2.2(c)]{Mig}.)

The following simple result due to Miller 
gives a connection between linkage theory and Alexander duality for monomial ideals defined in Definition \ref{alex}.

\begin{lemma}[{Miller \cite[Theorem 2.1]{Mi2}}]
\label{mlinkage}
Let $P=(x_1^{c_1+1},\dots,x_n^{c_n+1})$
and let $I \subset S$ be a $\cc$-ideal.
Then $P:(I+P)=I^\vee+P$.
\end{lemma}

\begin{lemma} \label{linkage}
For a $\cc$-ideal $I \subset S$, one has 
$$\pol(P):\pol(I+P)=\pol^*(I^\vee+P).$$
\end{lemma}

\begin{proof}
For any monomial $w \in \widetilde S$,
$w \in \pol(P):\pol(I+P)$ if and only if, for any $x^\aaa \in G(I+P)$,
$\pol(x^\aaa) w$ is divisible by $x_{i,0}x_{i,1}\cdots x_{i,c_i}$ for some $i$.
Then, since $\pol(x^\aaa)$ is a monomial of the form $ \prod_{a_i \ne 0} x_{i,0}x_{i,1} \cdots x_{i,a_i-1}$,
any generator $w \in G(\pol(P):\pol(I+P))$ must be of the form $w=\pol^*(x^\bb)$ for some $x^\bb \in S$.

On the other hand, for $\cc$-monomials $x^\aaa, x^\bb \in S$,
$\pol(x^\aaa)\pol^*(x^\bb) \in \pol(P)$ if and only if
$\pol(x^\aaa)\pol^*(x^\bb)$ is divisible by $x_{i,0}x_{i,1}\cdots x_{i,c_i}$ for some $i$,
equivalently, $x^\aaa x^\bb$ is divisible by $x_i^{c_i+1}$.
Thus $\pol(x^\aaa)\pol^*(x^\bb) \in \pol(P)$
if and only if $x^\aaa x^\bb \in P$.
Hence $\pol^*(x^\bb) \in G(\pol(P):\pol(I+P))$ if and only if $x^\bb \in G(P:(I+P))$.
Then the statement follows from Lemma \ref{mlinkage}. 
\end{proof}

Since $\pol(P)$
is generated by an $S$-sequence,
$\widetilde S/\pol(P)$ is a Gorenstein ring of dimension $|\cc|$.
Since the dimension of $S/(I+P)$ is $0$ and since taking polarizations preserves the Cohen-Macaulay property,
$\widetilde S/\pol(I+P)$ is a Cohen-Macaulay ring of dimension $|\cc|$.
Then the standard fact in linkage theory mentioned before Lemma \ref{mlinkage}
gives a purely algebraic proof of the following statement (apply the case when $Q=\pol(P)$ and $L=\pol(I+P)$).

\begin{theorem}
For any $\cc$-ideal $I \subset S$,
the ideal $\pol(I)+\pol^*(I^\vee)+\pol(P) \subset \widetilde S$ is Gorenstein.
\end{theorem}

\begin{remark}
By Lemmas \ref{complement} and \ref{linkage},
the Stanley-Reisner complex of $I_{\Lambda_\cc}:I_{\B_\cc(M)}=\pol(P): I_{\B_\cc(M)}$
is the complementary ball of $\B_\cc(M)$ in $\Lambda_\cc$.
This is the special case of the following general fact:
Let $\Gamma$ be a $d$-dimensional simplicial sphere on $[n]$,
$B \subset \Gamma$ a $d$-dimensional ball and
$B^c=\langle F \subset [n] : F \in \Gamma \setminus B \rangle$.
Then $I_\Gamma: I_B = I_{B^c}$.
\end{remark}

\subsection*{Graded Betti numbers and Bier spheres}
An important application of polarization appears in the study of graded Betti numbers.
Computations of the (multi) graded Betti numbers are one of current trends in combinatorial commutative algebra.
In the rest of this section, we study graded Betti numbers of Bier spheres.
We refer the readers to \cite{MS} for basics on multigraded commutative algebra.

Let $I \subset S$ be a monomial ideal.
Then $\Tor_i(S/I,K)$ is $\ZZ_{\geq 0}^n$-graded.
The integers
$$\beta_{i,(a_1,\dots,a_n)}=\dim_K \Tor_i(S/I,K)_{(a_1,\dots,a_n)}$$
are called the \textit{multigraded Betti numbers of $S/I$},
where $M_{(a_1,\dots,a_n)}$ denotes the homogeneous component of a $\ZZ_{\geq 0}^n$-graded $S$-module $M$ of degree $(a_1,\dots,a_n)$.
In the rest of this paper,
we identify $\aaa\in \ZZ_{\geq 0}^n$ with $x^\aaa$ for convenience.
The numbers $\beta_{i,j}(S/I)=\sum_{\deg x^\aaa=j} \beta_{i,x^\aaa}(S/I)$ and $\beta_i(S/I)=\sum_{j} \beta_{i,j}(S/I)$
are called \textit{graded Betti numbers} and \textit{total Betti numbers} of $S/I$ respectively.

The next result is useful to study graded Betti numbers of $K[\bier_\cc (M)]$.

\begin{lemma} \label{cone}
Let $B$ be a $(d-1)$-dimensional simplicial ball on $X=\{x_1,\dots,x_n\}$,
and let $\Delta =\partial B$ be the boundary complex of $B$.
If $B$ is a cone
then for any $i$ and for any squarefree monomial $x^F = \prod_{j \in F} x_j\in S$,
$$\beta_{i,x^F} (S/I_\Delta) = \beta_{i,x^F}(S/I_B) + \beta_{n+1-d-i, x^{[n] \setminus F}}(S/I_B).$$
In particular, $\beta_{i,j} (S/I_{\Delta}) = \beta_{i,j}(S/I_B) + \beta_{n+1-d-i, n-j}(S/I_B)$
for all $i,j$.
\end{lemma}

\begin{proof}
Since multigraded Betti numbers of Stanley-Reisner rings are concentrated in squarefree degrees (see \cite[Corollary 5.12]{MS}),
it is enough to prove the first statement.
The short exact sequence
$$0 \longrightarrow I_{\Delta}/I_B \longrightarrow S/I_B \longrightarrow S/I_{\Delta} \longrightarrow 0$$
yields the long exact sequence
\begin{align*}
&\cdots \longrightarrow \Tor_i(I_\Delta/I_B,K)_{x^F} 
\longrightarrow \Tor_i(S/I_B,K)_{x^F} 
\longrightarrow \Tor_i(S/I_\Delta,K)_{x^F} \\
&\hspace{18pt}
\longrightarrow \Tor_{i-1}(I_\Delta/I_B,K)_{x^F} 
\longrightarrow \Tor_{i-1}(S/I_B,K)_{x^F} 
\longrightarrow \cdots.
\end{align*}
If $B$ is a cone w.r.t. the vertex $x_k$ then all facets of $B$ contain $x_k$.
Then $G(I_B)$ has no monomials which are divisible by $x_k$.
Thus
\begin{align}
 \label{3.1}
\Tor_i(S/I_B,K)_{x^F}=0 \mbox{ if } k \in F.
\end{align}
On the other hand,
since $I_\Delta/I_B$ is the canonical module of $S/I_B$ \cite[II Theorem 7.3]{St},
we have
\begin{align} \label{3.2}
\Tor_i(I_\Delta/I_B,K)_{x^F}=\Tor_{n-d-i}(S/I_B,K)_{x^{[n]\setminus F}}
\end{align}
(see \cite[Theorem 13.37]{MS}) and
\begin{align} \label{3.3}
\Tor_i(I_\Delta/I_B,K)_{x^F}=0 \mbox{ if } k \not \in F.
\end{align}
Then by applying \eqref{3.1}, \eqref{3.2} and \eqref{3.3} to the long exact sequence we have
$$
\beta_{i,x^F}(S/I_\Delta)= \beta_{i,x^F}(S/I_B)
=\beta_{i,x^F}(S/I_B) + \beta_{n+1-d-i, x^{[n] \setminus F}}(S/I_B)
\ \mbox{ if } k \not \in F$$
and
$$
\beta_{i,x^F}(S/I_\Delta)=  \beta_{n+1-d-i, x^{[n] \setminus F}}(S/I_B)
=\beta_{i,x^F}(S/I_B) + \beta_{n+1-d-i, x^{[n] \setminus F}}(S/I_B)
\ \mbox{ if } k \in F,$$
as desired.
\end{proof}

Let $M$ be a $\cc$-multicomplex and $\mathrm{Lcm}(M)$ the least common multiple of monomials in $M$.
As we see in Remark \ref{r1}, if $\mathrm{Lcm}(M) \ne x^\cc$ then $\B_\cc(M)$ is a cone.
Since taking polarizations does not change graded Betti numbers,
Lemma \ref{cone} yields the next corollary.

\begin{corollary} \label{BettiNumber}
Let $M$ be a proper $\cc$-multicomplex on $X$.
If $x^\cc \ne \mathrm{Lcm}(M)$
then 
$$\beta^{\widetilde S}_{i,j}(\widetilde S/I_{\bier_\cc (M)}) = \beta_{i,j}^S(S/I(M)) + \beta^S_{n+1-i,|\overline \cc|-j}(S/I(M))$$
for all $i$ and $j$, where $\beta^{\widetilde S}_{i,j}(\widetilde S/I_{\bier_\cc (M)})$ are the graded Betti numbers over $\widetilde S$.
\end{corollary}

Note that the above formula does not always hold if $\mathrm{Lcm}(M) = x^\cc$.
For example, Example \ref{3..4} does not satisfy the formula.

Although we need an assumption on $\mathrm{Lcm}(M)$,
Corollary \ref{BettiNumber} will be useful to find many examples of graded Betti numbers of Gorenstein Stanley-Reisner rings.
For example, the above corollary implies the following non-trivial result.

\begin{corollary}
Let $I \subset S$ be a monomial ideal such that $S/I$ has finite length,
and let $b_i=\beta_i(S/I)$ for all $i$.
Then there exists a simplicial complex $\Delta$ such that
$K[\Delta]$ is Gorenstein and $\beta_i(K[\Delta])=b_i+b_{n+1-i}$ for all $i$.
\end{corollary}

\begin{proof}
If $S/I$ has finite length then there exists a finite multicomplex $M$ such that
$I=I(M)$. Choose a sufficiently large $\cc \in \ZZ_{\geq 0}^n$.
Corollary \ref{BettiNumber} says $\beta_i(K[\bier_\cc(M)])\!=b_i+b_{n+1-i}$ for all $i$.
\end{proof}

\begin{example} \label{BettiTable}
To understand the meaning of Corollary \ref{BettiNumber},
it is convenient to consider Betti tables (tables whose $i,j$-th entry is $\beta_{i,i+j}$).

Let $\cc=(2,2,2)$ and $M=\{1,x,y,z,xz,yz,z^2,yz^2\}$.
Then we have $I(M)=(x^2,y^2,z^3,xy,xz^2)$ and the Betti table of $K[x,y,z]/I(M)$
computed by the computer algebra system Macaulay 2 \cite{GS} is
\begin{verbatim}
     total: 1 5 6 2
         0: 1 . . .
         1: . 3 2 . 
         2: . 2 3 1
         3: . . 1 1
\end{verbatim}
Corollary \ref{BettiNumber} says that 
the Betti table of $K[\bier_\cc(M)]$ is the sum of the Betti table
of $K[x,y,z]/I(M)$ and the table which is obtained by transposing this table.
Thus the Betti table of $K[\bier_\cc(M)]$ is computed as follows:\medskip
\begin{verbatim}
     total: 1 5 6 2 .    total: . 2 6 5 1      total: 1 7 12 7 1
         0: 1 . . . .        0: . . . . .          0: 1 .  . . .
         1: . 3 2 . .        1: . . . . .          1: . 3  2 . .
         2: . 2 3 1 .  +     2: . 1 1 . .   =      2: . 3  4 1 .
         3: . . 1 1 .        3: . 1 3 2 .          3: . 1  4 3 .
         4: . . . . .        4: . . 2 3 .          4: . .  2 3 .
         5: . . . . .        5: . . . . 1          5: . .  . . 1
\end{verbatim}
\end{example}
\bigskip

\noindent
\textbf{Acknowledgments}:
This work was supported by KAKENHI 09J00756.
I would like to thank the referees for careful readings
and many helpful comments.

\end{document}